\documentclass[12pt,a4paper]{article}

\usepackage{theorem,enumerate}
\usepackage{amsmath,latexsym,amssymb,amsfonts}
\usepackage{eucal}
\usepackage{color}
\usepackage{comment}

\theorembodyfont{\normalfont\slshape}
\newtheorem{thm}{Theorem}

\newtheorem{prop}[thm]{Proposition}
\newtheorem{lem}[thm]{Lemma}

\newtheorem{claim}{Claim}[section]

\newtheorem{Thm}{Theorem}

\newcommand{\proof}{\medbreak\noindent\textit{Proof.}\quad}

\newcommand{\qed}{{$\quad\square$\vs{3.6}}}

\numberwithin{equation}{section}

\newcommand{\vs}[1]{\vspace*{#1 mm}}

\def\CC{{ \mathcal{C}}}

\def\GG{{ \mathcal{G}}}

\def\QQ{{ \mathcal{Q}}}

\bfseries\normalfont

\title{A new strategy for finding spanning trees without small degree stems}

\author{
Michitaka Furuya$^{1}$\footnote{\texttt{e-mail:mfuruya@kwansei.ac.jp}}\and \
Shoichi Tsuchiya$^{2}$\footnote{\texttt{e-mail:s.tsuchiya@isc.senshu-u.ac.jp}} \vs{5}\\
$^{1}$\textsl{School of Engineering,}\\
\textsl{Kwansei Gakuin University,}\\
\textsl{1 Gakuen Uegahara, Sanda, Hyogo, 669-1330, Japan}\\
$^{3}$\textsl{School of Network and Information,}\\
\textsl{Senshu University,}\\
\textsl{2-1-1 Higashimita, Tama-ku, Kawasaki-shi, Kanagawa, 214-8580, Japan}\\
}

\date{}

\begin{document}

\maketitle

\begin{abstract}
For an integer $k\geq 2$, a spanning tree of a graph without vertices of degree from $2$ to $k$ is called a {\it $[2,k]$-ST} of the graph.
The concept of $[2,k]$-STs is a natural extension of a homeomorphically irreducible spanning tree (or HIST), which is a well-studied graph structure.
In this paper, we give a new strategy for finding $[2,k]$-STs.
By using the strategy, we refine or extend a known degree-sum condition for the existence of a HIST.
Furthermore, we also investigate a degree-product condition for the existence of a $[2,k]$-ST.
\end{abstract}



\section{Introduction}\label{sec1}

For a graph $G$, let $V(G)$ and $E(G)$ denote the {\it vertex set} and the {\it edge set} of $G$, respectively.
For $u\in V(G)$, let $N_{G}(u)$ and $d_{G}(u)$ denote the {\it neighborhood} and the {\it degree} of $u$, respectively; thus $N_{G}(u)=\{v\in V(G):uv\in E(G)\}$ and $d_{G}(u)=|N_{G}(u)|$.
For an integer $i\geq 0$, let $V_{i}(G)=\{u\in V(G):d_{G}(u)=i\}$ and $V_{\geq i}(G)=\{u\in V(G):d_{G}(u)\geq i\}$.
We let $\delta (G)$ denote the {\it minimum degree} of $G$.
We define 
$$
\sigma _{2}(G)=\min\{d_{G}(u)+d_{G}(v):u,v\in V(G),~u\neq v,~uv\notin E(G)\}
$$
if $G$ is not complete; 
$\sigma _{2}(G)=\infty $ if $G$ is complete.
Also, we define
$$
\pi _{2}(G)=\min\{d_{G}(u)d_{G}(v):u,v\in V(G),~u\neq v,~uv\notin E(G)\}
$$
if $G$ is not complete; 
$\pi _{2}(G)=\infty $ if $G$ is complete.

 For a tree $T$, 
each vertex in $V_{\ge 2}$ (resp., $V_{1}(T)$) is called a {\it stem} (resp., a {\it leaf}). 
For a graph $G$, 
a spanning tree of $G$ without vertices of degree $2$ is called a {\it homeomorphically irreducible spanning tree} (or a {\it HIST}) of $G$; i.e., a spanning tree $T$ of $G$ is a HIST if and only if $V_{2}(T)=\emptyset $.
A structure of HISTs is sometimes used as an essential tool to construct graph classes; for example, in an explicit class of edge-minimal $3$-connected plane graphs given by Halin~\cite{H}, HISTs play a key role.
Motivated by such uses, 
the existence of a HIST (or a large subtree having no vertex of degree $2$) has been widely studied (for example, see \cite{ABHT,CRS,CS,HNO,ZWHY}).
It is well known that a number of sufficient conditions for the existence
of a hamiltonian path have been naturally generalized to those for the existence
of a {\it spanning $k$-tree}, which is a spanning tree in which every stem has degree lying between $2$ and $k$.   
Similar to this, the concept of HISTs was naturally extended:
A spanning tree $T$ of $G$ is called a {\it $[2,k]$-ST} of $G$ if $\bigcup _{2\leq i\leq k}V_{i}(T)=\emptyset $ (for further historical background and related results, we refer the reader to \cite{FST}).

Our aim in this series is to refine and to extend some known degree conditions for the existence of HISTs.
In this paper,
\begin{enumerate}[{\upshape(i)}]
\item
we give two results which essentially extend a known degree-sum condition assuring us the existence of a HIST, and
\item
we focus on a degree-product condition, which seems to be more reasonable for the existence of a HIST, and find a $[2,k]$-ST using such a condition.
\end{enumerate}

We start with a degree-sum condition for the existence of HISTs, which was recently given by Ito and Tsuchiya~\cite{IT}.

\begin{Thm}[Ito and Tsuchiya~\cite{IT}]
\label{ThmB}
Let $G$ be a connected graph of order $n\geq 8$.
If $\sigma _{2}(G)\geq n-1$, then $G$ has a HIST.
\end{Thm}

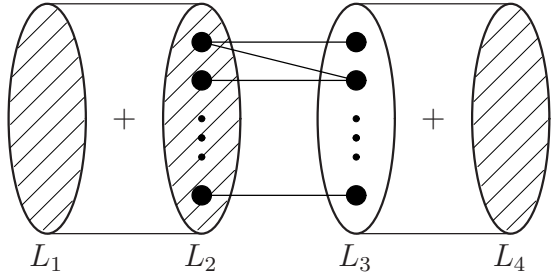
\begin{figure}
\begin{center}
{\unitlength 0.1in%
\begin{picture}(28.0000,12.6000)(4.0000,-14.6000)%
%
\special{pn 13}%
\special{ar 600 800 200 600 0.0000000 6.2831853}%
%
\special{pn 13}%
\special{ar 1400 800 200 600 0.0000000 6.2831853}%
%
\special{pn 13}%
\special{ar 2200 800 200 600 0.0000000 6.2831853}%
%
\special{pn 13}%
\special{ar 3000 800 200 600 0.0000000 6.2831853}%
%
\special{pn 8}%
\special{pa 2200 200}%
\special{pa 3000 200}%
\special{fp}%
%
\special{pn 8}%
\special{pa 2200 1400}%
\special{pa 3000 1400}%
\special{fp}%
%
\special{pn 8}%
\special{pa 600 1400}%
\special{pa 1400 1400}%
\special{fp}%
%
\special{pn 8}%
\special{pa 600 200}%
\special{pa 1400 200}%
\special{fp}%
\put(5.9500,-15.2500){\makebox(0,0){$L_{1}$}}%
\put(10.0000,-8.0000){\makebox(0,0){$+$}}%
\put(26.0000,-8.0000){\makebox(0,0){$+$}}%
\put(13.9500,-15.2500){\makebox(0,0){$L_{2}$}}%
\put(21.9500,-15.2500){\makebox(0,0){$L_{3}$}}%
\put(29.9500,-15.2500){\makebox(0,0){$L_{4}$}}%
%
\special{sh 1.000}%
\special{ia 1400 400 50 50 0.0000000 6.2831853}%
\special{pn 8}%
\special{ar 1400 400 50 50 0.0000000 6.2831853}%
%
\special{sh 1.000}%
\special{ia 1400 600 50 50 0.0000000 6.2831853}%
\special{pn 8}%
\special{ar 1400 600 50 50 0.0000000 6.2831853}%
%
\special{sh 1.000}%
\special{ia 1400 1200 50 50 0.0000000 6.2831853}%
\special{pn 8}%
\special{ar 1400 1200 50 50 0.0000000 6.2831853}%
%
\special{pn 4}%
\special{sh 1}%
\special{ar 1400 800 16 16 0 6.2831853}%
\special{sh 1}%
\special{ar 1400 1000 16 16 0 6.2831853}%
\special{sh 1}%
\special{ar 1400 900 16 16 0 6.2831853}%
\special{sh 1}%
\special{ar 1400 900 16 16 0 6.2831853}%
%
\special{sh 1.000}%
\special{ia 2200 400 50 50 0.0000000 6.2831853}%
\special{pn 8}%
\special{ar 2200 400 50 50 0.0000000 6.2831853}%
%
\special{sh 1.000}%
\special{ia 2200 600 50 50 0.0000000 6.2831853}%
\special{pn 8}%
\special{ar 2200 600 50 50 0.0000000 6.2831853}%
%
\special{sh 1.000}%
\special{ia 2200 1200 50 50 0.0000000 6.2831853}%
\special{pn 8}%
\special{ar 2200 1200 50 50 0.0000000 6.2831853}%
%
\special{pn 4}%
\special{sh 1}%
\special{ar 2200 800 16 16 0 6.2831853}%
\special{sh 1}%
\special{ar 2200 1000 16 16 0 6.2831853}%
\special{sh 1}%
\special{ar 2200 900 16 16 0 6.2831853}%
\special{sh 1}%
\special{ar 2200 900 16 16 0 6.2831853}%
%
\special{pn 8}%
\special{pa 1400 400}%
\special{pa 2200 400}%
\special{fp}%
%
\special{pn 8}%
\special{pa 1400 400}%
\special{pa 2200 600}%
\special{fp}%
%
\special{pn 8}%
\special{pa 1400 1200}%
\special{pa 2200 1200}%
\special{fp}%
%
\special{pn 8}%
\special{pa 1400 600}%
\special{pa 2200 600}%
\special{fp}%
%
\special{pn 4}%
\special{pa 790 690}%
\special{pa 430 1050}%
\special{fp}%
\special{pa 790 810}%
\special{pa 450 1150}%
\special{fp}%
\special{pa 790 930}%
\special{pa 470 1250}%
\special{fp}%
\special{pa 770 1070}%
\special{pa 510 1330}%
\special{fp}%
\special{pa 730 1230}%
\special{pa 580 1380}%
\special{fp}%
\special{pa 780 580}%
\special{pa 410 950}%
\special{fp}%
\special{pa 760 480}%
\special{pa 410 830}%
\special{fp}%
\special{pa 740 380}%
\special{pa 410 710}%
\special{fp}%
\special{pa 700 300}%
\special{pa 420 580}%
\special{fp}%
\special{pa 650 230}%
\special{pa 450 430}%
\special{fp}%
%
\special{pn 4}%
\special{pa 1590 690}%
\special{pa 1230 1050}%
\special{fp}%
\special{pa 1590 810}%
\special{pa 1250 1150}%
\special{fp}%
\special{pa 1590 930}%
\special{pa 1270 1250}%
\special{fp}%
\special{pa 1570 1070}%
\special{pa 1310 1330}%
\special{fp}%
\special{pa 1530 1230}%
\special{pa 1380 1380}%
\special{fp}%
\special{pa 1580 580}%
\special{pa 1210 950}%
\special{fp}%
\special{pa 1560 480}%
\special{pa 1210 830}%
\special{fp}%
\special{pa 1540 380}%
\special{pa 1210 710}%
\special{fp}%
\special{pa 1500 300}%
\special{pa 1220 580}%
\special{fp}%
\special{pa 1450 230}%
\special{pa 1250 430}%
\special{fp}%
%
\special{pn 4}%
\special{pa 3190 690}%
\special{pa 2830 1050}%
\special{fp}%
\special{pa 3190 810}%
\special{pa 2850 1150}%
\special{fp}%
\special{pa 3190 930}%
\special{pa 2870 1250}%
\special{fp}%
\special{pa 3170 1070}%
\special{pa 2910 1330}%
\special{fp}%
\special{pa 3130 1230}%
\special{pa 2980 1380}%
\special{fp}%
\special{pa 3180 580}%
\special{pa 2810 950}%
\special{fp}%
\special{pa 3160 480}%
\special{pa 2810 830}%
\special{fp}%
\special{pa 3140 380}%
\special{pa 2810 710}%
\special{fp}%
\special{pa 3100 300}%
\special{pa 2820 580}%
\special{fp}%
\special{pa 3050 230}%
\special{pa 2850 430}%
\special{fp}%
\end{picture}}%
\caption{Graphs $G$ belonging to $\GG_{k,n}$}
\label{f2}
\end{center}
\end{figure}

They also showed that the bound on $\sigma _{2}$ is best possible, i.e., for each integer $n\geq 8$, there exists a graph $G$ of order $n$ with $\sigma _{2}(G)=n-2$ having no HIST.
Our first result is a refinement of Theorem~\ref{ThmB} with a characterization of sharp examples.
For integers $k\geq 2$ and $n\geq 2k+1$, let $\GG_{k,n}$ be the family of graphs $G$ of order $n$ satisfying the following conditions (see Figure~\ref{f2}):
\begin{enumerate}[{\bf (L1)}]
\setlength{\parskip}{0cm}
\setlength{\itemsep}{0cm}
\item[{\bf (L1)}]
$V(G)$ is the disjoint union of four non-empty sets $L_{1}$, $L_{2}$, $L_{3}$ and $L_{4}$,
\item[{\bf (L2)}]
$L_{1}\cup L_{2}$ and $L_{4}$ are cliques of $G$,
\item[{\bf (L3)}]
for every $u_{1}\in L_{1}$, $N_{G}(u_{1})\cap (L_{3}\cup L_{4})=\emptyset $,
\item[{\bf (L4)}]
for every $u_{2}\in L_{2}$, $N_{G}(u_{2})\cap L_{3}\neq \emptyset $, $N_{G}(u_{2})\cap L_{4}=\emptyset $ and $d_{G}(u_{2})\leq k$, and
\item[{\bf (L5)}]
for every $u_{3}\in L_{3}$, $N_{G}(u_{3})\cap L_{2}\neq \emptyset $, $L_{4}\subseteq N_{G}(u_{3})$ and $d_{G}(u_{3})\geq n-|L_{1}\cup L_{2}|-1$.
\end{enumerate}
Let $k\geq 2$ be an integer, and let $c_{k}=\sqrt{k(k-1)(k+2\sqrt{2k}+2)}$.
Let $n_{0}(k)$ be the smallest positive integer such that $n-4c_{k}\sqrt{n}-2k^{2}-4k-4\geq 0$ for every integer $n\geq n_{0}(k)$.
Our first result is the following.

\begin{thm}
\label{thm2+}
Let $k\geq 2$ be an integer.
Let $G$ be a connected graph of order $n\geq n_{0}(k)$, and suppose that $\sigma _{2}(G)\geq n-2$.
Then $G$ has a $[2,k]$-ST if and only if $G$ is not isomorphic to any graph in $\GG_{k,n}$.
\end{thm}

Theorem~\ref{thm2+} is a generalization of Theorem~\ref{ThmB} for sufficiently large graphs.
On the other hand, one can easily calculate that $n_{0}(2)=295$, and so it in fact does not cover Theorem~\ref{ThmB} when a target graph is small.
Actually, in the previous paper~\cite{FST} of this series, we obtained the same result as Theorem~\ref{thm2+} for the case where $k=2$ and $n\geq 10$ in a different way.
So the order condition $n\geq n_{0}(k)$ in Theorem~\ref{thm2+} is not best possible.

Recently, Shan and Tsuchiya~\cite{ST} introduced a blocking set, which is a new concept on cutsets closely related to the existence of a HIST.
We extend the concept to a $[2,k]$-ST version.
Let $k\geq 2$ be an integer, and suppose that $G$ is connected.
A cutset $U\subseteq V(G)$ of $G$ is {\it $k$-blocking set} of $G$ if $U\subseteq \bigcup _{2\leq i\leq k}V_{i}(G)$.
If a graph $G$ has a $[2,k]$-ST $T$, then for a cutset $L$ of $G$, there exists a vertex $u\in L$ with $d_{T}(u)\geq k+1$.
In particular, if a graph has a $[2,k]$-ST, then the graph has no $k$-blocking set.

If a graph $G$ satisfies (L1)--(L5), then $L_{2}$ is a $k$-blocking set.
Thus, considering Theorem~\ref{thm2+}, one might expect that the degree-sum condition can be greatly improved if we omit the existence of a $k$-blocking set.
Our second result affirms the expectation.
Let $n_{1}(k)$ be the smallest positive integer such that $\frac{n+2k-2}{4}-2c_{k}\sqrt{n}-k^{2}-2k-1\geq 0$ for every integer $n\geq n_{1}(k)$.
Note that $n_{1}(2)=1091$.

\begin{thm}
\label{thm3}
Let $k\geq 2$ be an integer.
Let $G$ be a connected graph of order $n\geq n_{1}(k)$, and suppose that $\sigma _{2}(G)\geq \frac{n+2k-2}{2}$.
Then $G$ has a $[2,k]$-ST if and only if $G$ has no $k$-blocking set.
\end{thm}

Our third result is to propose a new concept on degree conditions.
To explain it in detail, we start with two more natural results on degree conditions for the existence of HISTs (or $[2,k]$-STs).
The following theorem is the first result discussing a relationship between a HIST and a degree condition.

\begin{Thm}[Albertson, Berman, Hutchinson and Thomassen~\cite{ABHT}]
\label{ThmA}
Let $G$ be a connected graph of order $n$, and suppose that $\delta (G)\geq 4\sqrt{2n}$.
Then $G$ has a HIST.
\end{Thm}

Note that $c_{2}=4$ because $c_{k}=\sqrt{k(k-1)(k+2\sqrt{2k}+2)}$.
Recently, Theorem~\ref{ThmA} was refined and extended in the previous paper of this series as follows.

\begin{Thm}[Furuya, Saito and Tsuchiya~\cite{FST}]
\label{ThmC}
Let $k\geq 2$ be an integer.
Let $G$ be a connected graph of order $n$, and suppose that $\delta (G)\geq c_{k}\sqrt{n}$.
Then $G$ has a $[2,k]$-ST.
\end{Thm}

The coefficient of $\sqrt{n}$ in Theorem~\ref{ThmC} might be further improved.
On the other hand, for any integers $k\geq 2$ and $d\geq k-1$ such that $\frac{d}{k-1}$ is an integer, Furuya et al.~\cite{FST} constructed a connected graph $G$ with $\delta (G)=d=\sqrt{4(k-1)|V(G)|+(2k-1)^{2}}-2k+1$ having no $[2,k]$-ST.
Therefore the degree condition in Theorem~\ref{ThmC} is asymptotically best possible.

Now we focus on a large gap between Theorems~\ref{ThmB} and \ref{ThmA}.
For example, the following two theorems are well-known, and their degree conditions are best possible:
\begin{enumerate}[{$\bullet $}]
\item
Dirac's Theorem~\cite{Di}: If a graph $G$ of order $n\geq 3$ satisfies $\delta (G)\geq \frac{n}{2}$, then $G$ has a Hamiltonian cycle.
\item
Ore's Theorem~\cite{O}: If a graph $G$ of order $n\geq 3$ satisfies $\sigma _{2}(G)\geq n$, then $G$ has a Hamiltonian cycle.
\end{enumerate}
In particular, Ore's Theorem implies Dirac's Theorem.
On the other hand, a degree condition in Theorem~\ref{ThmB} is much bigger than one in Theorem~\ref{ThmA}.
Considering the fact that the root of the order of a graph appears in Theorem~\ref{ThmA}, one natural question occurs:
Is there a degree-product condition close to the order assuring us the existence of HISTs?
We give an affirmative answer for the problem.
Recall that $c_{k}=\sqrt{k(k-1)(k+2\sqrt{2k}+2)}$ for an integer $k\geq 2$.
For an integer $k\geq 2$, let
$$
p_{k}=\frac{5c_{k}^{2}+3c_{k}\sqrt{c_{k}^{2}+4k^{2}+8k+4}}{2}+k^{2}+2k+1.
$$
Our third result is the following.

\begin{thm}
\label{thm4}
Let $k\geq 2$ be an integer.
Let $G$ be a connected graph of order $n\geq k+2$, and suppose that $\pi _{2}(G)\geq p_{k}n$.
Then $G$ has a $[2,k]$-ST.
\end{thm}

As with Theorems~\ref{ThmA} and \ref{ThmC}, we do not know whether the coefficient of $n$ in Theorem~\ref{thm4} is best possible or not.
However, the degree-product condition, which is an unprecedented work as we know, seems to be essential for the existence of HISTs (or $[2,k]$-STs).

The proofs of Theorem~\ref{thm2+}--\ref{thm4} depend on a common strategy.
In Section~\ref{sec-keylem}, we introduce key lemmas for the strategy.
In Section~\ref{sec-pf-thm3}, we prove Theorems~\ref{thm2+} and \ref{thm3} at the same time, and discuss the sharpness of the degree-sum condition appearing in Theorem~\ref{thm3}.
In Section~\ref{sec-pf-thm4}, we prove Theorem~\ref{thm4}.

\subsection{Notations}\label{sec-pre}

In this subsection, we prepare the notation required for our proofs.
For terms and symbols not defined in this paper, we refer the reader to \cite{D}.

Let $G$ be a graph.
For $F\subseteq E(G)$, let $V(F)=\{u,v:uv\in F\}$.
For a subgraph $H$ of $G$ and a subset $F$ of $E(G)$, let $H+F$ be the subgraph of $G$ with $V(H+F)=V(H)\cup V(F)$ and $E(H+F)=E(H)\cup F$.
Let ${\rm compo}(G)$ be the number of components of $G$.
A vertex $u$ of $G$ is called a {\it cut-vertex} of $G$ if ${\rm compo} (G-u)>{\rm compo} (G)$.
Note that cut-vertices are defined in disconnected graphs.
Let ${\rm cut}(G)$ be the number of cut-vertices of $G$.

Let $k\geq 2$ be an integer.
For a tree $T$ and a subset $U$ of $V(T)$, $T$ is {\it $(k,U)$-good} if $V(T)\setminus U\subseteq V_{1}(T)\cup V_{\geq k+1}(T)$ and $U\subseteq V_{\geq k}(T)$.
Note that a spanning $(k,\emptyset )$-good tree of a graph $G$ is a $[2,k]$-ST of $G$.
If a tree is $(k,\{u\})$-good, then the tree is simply said to be {\it $(k,u)$-good}.

\section{Key lemmas}\label{sec-keylem}

In this section, we introduce a key lemma for our argument (Lemma~\ref{lem-3-1++}) and arrange it for the existence of $[2,k]$-STs (Lemmas~\ref{lem-3-2} and \ref{lem-3-3}).
Our strategy is that first we take the vertex set $S$ in a graph $G$ consisting of all small degree vertices, where small means half or the root of the degree condition (by the definitions of degree conditions, we can show that $S$ induces a clique). Then we can see that each component of $G-S$ has large minimum degree. In order to take $[2,k]$-STs of $G$, 
we guarantee the existence of  convenient structures in such components by proving Lemmas~\ref{lem-3-2} and \ref{lem-3-3}. 

\begin{lem}
\label{lem-3-1}
Let $G$ be a graph of order $n$, and suppose that $\delta (G)\geq 2\sqrt{n}$.
Then ${\rm cut}(G)+{\rm compo} (G)-1\leq 2\sqrt{n}$.
\end{lem}
\proof
We proceed by induction on $n$.
If $n\leq 5$, then there is no graph $G$ of order $n$ with $\delta (G)\geq 2\sqrt{n}$, and hence the lemma holds.
Thus we may assume that $n\geq 6$.

Since $\delta (G)\geq 2\sqrt{n}$, every component of $G$ contains more than $2\sqrt{n}$ vertices, and hence ${\rm compo} (G)<\frac{n}{2\sqrt{n}}=\frac{\sqrt{n}}{2}$.
In particular, if $G$ has no cut-vertex, then ${\rm cut}(G)+{\rm compo} (G)-1<0+\frac{\sqrt{n}}{2}-1<2\sqrt{n}$, as desired.
Thus we may assume that a component $G_{1}$ of $G$ has a cut-vertex.

Let $L$ be an {\it end-block} of $G_{1}$, which is a block of $G$ containing exactly one cut-vertex.  
Let $u$ be the unique cut-vertex of $G_{1}$ with $u\in V(L)$.
For a vertex $u'\in V(L)\setminus \{u\}$, $2\sqrt{n}\leq \delta (G)\leq d_{G}(u')\leq |V(L)\setminus \{u'\}|=|V(L)|-1$.
Hence
\begin{align}
|V(L)|\geq 2\sqrt{n}+1.\label{cond-lem-3-1-1}
\end{align}
Since $u$ is a cut-vertex of $G_{1}$, $V(G_{1})\setminus V(L)\neq \emptyset $.
Furthermore, since $L$ is an end-block, all cut-vertices of $G_{1}$ other than $u$ are contained in $V(G_{1})\setminus V(L)$.
Let $X$ be the set of cut-vertices of $G_{1}$ other than $u$ such that they are not cut-vertices of $G_{1}-V(L)$.

Fix a vertex $v\in X$, and let $H$ be the component of $G_{1}-V(L)$ containing $v$.
Then $v$ belongs to exactly two blocks of $G_{1}$ and all neighbors of $v$ in one of them have been deleted in $G-V(L)$.
This implies that $uv\in E(G)$ and $G_{1}[\{u,v\}]$ is a block of $G_{1}$.
In particular, $V(H)\cap X=\{v\}$.
Since $v$ is arbitrary, $|X|=|\{H:H\mbox{ is the component of }G_{1}-V(L)\mbox{ containing a vertex in }X\}|\leq {\rm compo} (G_{1}-V(L))$.
This implies that
\begin{align}
{\rm cut}&(G-V(L))+{\rm compo} (G-V(L))+1\nonumber \\
&= ({\rm cut}(G-V(G_{1}))+{\rm cut}(G_{1}-V(L)))+({\rm compo} (G)-|\{G_{1}\}|+{\rm compo} (G_{1}-V(L)))+1\nonumber \\
&\geq {\rm cut}(G-V(G_{1}))+({\rm cut}(G_{1})-|\{u\}\cup X|)+{\rm compo} (G)-1+|X|+1\nonumber \\
&= {\rm cut}(G)+{\rm compo} (G)-1.\label{cond-lem-3-1-2}
\end{align}
Since $(2\sqrt{n}-2)^{2}-\left(2\sqrt{n-2\sqrt{n}-1}\right)^{2}=8>0$, we have
\begin{align}
2\sqrt{n}-2>2\sqrt{n-2\sqrt{n}-1}.\label{cond-lem-3-1-3}
\end{align}
By (\ref{cond-lem-3-1-1}) and (\ref{cond-lem-3-1-3}), $\delta (G-V(L))\geq \delta (G)-1\geq 2\sqrt{n}-1>2\sqrt{n-2\sqrt{n}-1}\geq 2\sqrt{|V(G-V(L))|}$.
Hence, by the induction hypothesis on $G-V(L)$, (\ref{cond-lem-3-1-1}) and (\ref{cond-lem-3-1-3}), we have
$$
{\rm cut}(G-V(L))+{\rm compo} (G-V(L))\leq 2\sqrt{|V(G-V(L))|}+1\leq 2\sqrt{n-2\sqrt{n}-1}+1<2\sqrt{n}-1,
$$
This together with (\ref{cond-lem-3-1-2}) leads to
$$
{\rm cut}(G)+{\rm compo} (G)-1\leq {\rm cut}(G-V(L))+{\rm compo} (G-V(L))+1<(2\sqrt{n}-1)+1,
$$
as desired.
\qed

\begin{lem}
\label{lem-3-1++}
Let $m\geq 0$ be an integer.
Let $G$ be a connected graph of order $n$, and let $u\in V(G)$ and $Y\subseteq V(G)\setminus \{u\}$.
If $\delta (G)\geq 2\sqrt{n}+m+|Y|$, then there exists a set $X\subseteq N_{G}(u)\setminus Y$ with $|X|=m$ such that $G-X$ is connected.
\end{lem}
\proof
We proceed by induction on $m$.
If $m=0$, then the desired conclusion clearly holds.
Thus we may assume that $m\geq 1$.

Since $G$ is connected, it follows from Lemma~\ref{lem-3-1} that ${\rm cut}(G)\leq 2\sqrt{n}$.
Since $d_{G}(u)\geq 2\sqrt{n}+m+|Y|$, this implies that there exists a vertex $v\in N_{G}(u)\setminus Y$ which is not a cut-vertex of $G$.
Let $G'=G-v$.
Then $G'$ is connected and $\delta (G')\geq \delta (G)-1\geq 2\sqrt{n}+(m-1)+|Y|>2\sqrt{|V(G')|}+(m-1)+|Y|$.
Hence by the induction hypothesis on $G'$, there exists a set $X'\subseteq N_{G'}(u)\setminus Y~(=N_{G}(u)\setminus (Y\cup \{v\}))$ with $|X'|=m-1$ such that $G'-X'~(=G-(\{v\}\cup X'))$ is connected.
Consequently, $X:=\{v\}\cup X'$ is a desired subset of $N_{G}(u)\setminus Y$.
\qed

In the remainder of this section, we implicitly use the fact that $c_{k}>2$ for every integer $k\geq 2$.

\begin{lem}
\label{lem-3-2}
Let $k\geq 2$ be an integer.
Let $G$ be a connected graph of order $n$, and let $U\subseteq V(G)$ be a set with $U\neq \emptyset $.
If $\delta (G)\geq c_{k}\sqrt{n}+(k+1)|U|-1$, then there exists a spanning forest of $G$ consisting of exactly $|U|$ components $F_{1},F_{2},\ldots ,F_{|U|}$ such that for every integer $i$ with $1\leq i\leq |U|$, $|V(F_{i})\cap U|=1$ and $F_{i}$ is a $(k,V(F_{i})\cap U)$-good tree.
\end{lem}
\proof
Write $U=\{u_{1},u_{2},\ldots ,u_{t}\}$ where $t=|U|$.
Since $(k+1)t-(k+t)=k(t-1)\geq 0$, we have $\delta (G)\geq c_{k}\sqrt{n}+(k+1)t-1>2\sqrt{n}+k+t-1=2\sqrt{n}+k+|U\setminus \{u_{t}\}|$.
This together with Lemma~\ref{lem-3-1++} with $(m,u,Y)=(k,u_{t},U\setminus \{u_{t}\})$ implies that there exists a set $X\subseteq N_{G}(u_{t})\setminus (U\setminus \{u_{t}\})~(=N_{G}(u_{t})\setminus U)$ with $|X|=k$ such that $G-X$ is connected.

We proceed by induction on $t$.
Suppose that $t=1$, i.e., $U=\{u_{1}\}$.
Then $\delta (G-X)\geq \delta (G)-k\geq c_{k}\sqrt{n}>c_{k}\sqrt{|V(G-X)|}$.
Hence by Theorem~\ref{ThmC}, $G-X$ has a $[2,k]$-ST $T_{0}$.
Since $X\subseteq N_{G}(u_{1})$, $F_{1}:=T_{0}+\{u_{1}v:v\in X\}$ is a $[2,k]$-ST of $G$, and in particular, $F_{1}$ is a spanning $(k,u_{1})$-good tree of $G$, and hence it is a desired forest.
Thus we may assume that $t\geq 2$.

Let $\CC$ be the family of components of $G-(\{u_{t}\}\cup X)$, and let $\CC_{1}=\{C\in \CC:V(C)\cap U\neq \emptyset \}$ and $\CC_{2}=\CC\setminus \CC_{1}$.
Since $t\geq 2$, we have $\CC_{1}\neq \emptyset $.

Fix $C\in \CC_{1}$.
Let $I_{C}=\{i:1\leq i\leq t-1,~u_{i}\in V(C)\}$, and let $t_{C}=|I_{C}|$.
Then
\begin{align*}
\delta (C) &\geq \delta (G)-|\{u_{t}\}\cup X|\\
&\geq c_{k}\sqrt{n}+(k+1)t-1-(k+1)\\
&= c_{k}\sqrt{n}+(k+1)(t-1)-1\\
&>c_{k}\sqrt{|V(C)|}+(k+1)t_{C}-1.
\end{align*}
By the induction hypothesis on $C$, $C$ has a spanning forest consisting of exactly $t_{C}$ components $F_{i}~(i\in I_{C})$ such that for every integer $i\in I_{C}$, $|V(F_{i})\cap U|=1$ and $F_{i}$ is a $(k,V(F_{i})\cap U)$-good tree.
Since $C$ is arbitrary, $\sum _{C\in \CC_{1}}t_{C}~(=t-1)$ vertex-disjoint subtrees $F_{1},F_{2},\ldots ,F_{t-1}$ of $G$ have been defined.
Note that $\bigcup _{1\leq i\leq t-1}V(F_{i})=\bigcup _{C\in \CC_{1}}V(C)=V(G)\setminus (\{u_{t}\}\cup X\cup (\bigcup _{C'\in \CC_{2}}V(C')))$.

Remark that $\CC_{2}$ might be empty.
Assume that $\CC_{2}\neq \emptyset $ and fix $C'\in \CC_{2}$.
Since $G-X$ is connected, there exists a vertex $v_{C'}\in N_{G}(u_{t})\cap V(C')$.
Since $(k+1)t-1-(k+1)-k=(k+1)(t-2)\geq 0$, $(k+1)t-1-(k+1)\geq k$, and hence
$$
\delta (C')\geq \delta (G)-|\{u_{t}\}\cup X|\geq c_{k}\sqrt{n}+(k+1)t-1-(k+1)>c_{k}\sqrt{|V(C')|}+k.
$$
By the induction hypothesis on $C'$, $C'$ has a spanning $(k,v_{C'})$-good tree $T_{C'}$.
Let $F'$ be the subgraph of $G$ with $V(F')=\{u_{t}\}\cup X$ and $E(F')=\{u_{t}v:v\in X\}$.
Then $F'$ is a $(k,u_{t})$-good subtree of $G$.
Let $F_{t}=(F'\cup (\bigcup _{C'\in \CC_{2}}T_{C'}))+\{u_{t}v_{C'}:C'\in \CC_{2}\}$, where $F_{t}=F'$ if $\CC_{2}=\emptyset $.
Then $F_{t}$ is a $(k,u_{t})$-good tree, $V(F_{t})=\{u_{t}\}\cup X\cup (\bigcup _{C'\in \CC_{2}}V(C'))$ and $V(F_{t})\cap U=\{u_{t}\}$.

Consequently, the graph $\bigcup _{1\leq i\leq t}F_{i}$ is a desired spanning forest of $G$.
\qed

\begin{lem}
\label{lem-3-3}
Let $k\geq 2$ be an integer.
Let $G$ be a connected graph of order $n$, and let $U\subseteq V(G)$ be a set with $U\neq \emptyset $.
If $\delta (G)\geq c_{k}\sqrt{n}+k|U|-1$, then there exists a spanning $(k,U)$-good tree of $G$.
\end{lem}
\proof
Write $U=\{u_{1},u_{2},\ldots ,u_{t}\}$ where $t=|U|$.
We recursively define $t$ sets $X_{1},X_{2},\ldots ,X_{t}$ such that for each integer $i$ with $1\leq i\leq t$,
\begin{enumerate}[{\bf (A1)}]
\item[{\bf (A1)}]
$X_{i}\subseteq N_{G}(u_{i})\setminus (U\cup (\bigcup _{1\leq j\leq i-1}X_{j}))$,
\item[{\bf (A2)}]
$|X_{i}|=k-1$, and
\item[{\bf (A3)}]
$G-(\bigcup _{1\leq j\leq i}X_{j})$ is connected
\end{enumerate}
as follows:
We let $i_{0}$ be an integer with $1\leq i_{0}\leq t$, and assume that we have defined $i_{0}-1$ sets $X_{1},X_{2},\ldots ,X_{i_{0}-1}$ satisfying (A1)--(A3) for every integer $i$ with $1\leq i\leq i_{0}-1$.
Let $G_{i_{0}}=G-(\bigcup _{1\leq i\leq i_{0}-1}X_{i})$.
Then $G_{i_{0}}$ is connected and
\begin{align*}
\delta (G_{i_{0}}) &\geq \delta (G)-\left|\sum _{1\leq j\leq i_{0}-1}X_{j}\right|\\
&\geq c_{k}\sqrt{n}+kt-1-(i_{0}-1)(k-1)\\
&> 2\sqrt{n}+kt-1-(t-1)(k-1)\\
&= 2\sqrt{n}+k-1+|U\setminus \{u_{i_{0}}\}|.
\end{align*}
This together with Lemma~\ref{lem-3-1++} with $(G,m,u,Y)=(G_{i_{0}},k-1,u_{i_{0}},U\setminus \{u_{i_{0}}\})$ implies that there exists a set $X_{i_{0}}\subseteq N_{G_{i_{0}}}(u_{i_{0}})\setminus (U\setminus \{u_{i_{0}}\})~(=N_{G_{i_{0}}}(u_{i_{0}})\setminus U)$ with $|X_{i_{0}}|=k-1$ such that $G_{i_{0}}-X_{i_{0}}~(=G-(\bigcup _{1\leq i\leq i_{0}}X_{i}))$ is connected.
Thus $X_{1},X_{2},\ldots ,X_{i_{0}}$ satisfy (A1)--(A3) for for every integer $i$ with $1\leq i\leq i_{0}$.
Consequently, we obtain desired sets.

Let $G'=G-(\bigcup _{1\leq i\leq t}X_{i})$.
By (A3), $G'$ is connected.
Since $X_{1},X_{2},\ldots ,X_{t}$ are pairwise disjoint by (A1), it follows from (A2) that $\delta (G')\geq \delta (G)-|\sum _{1\leq i\leq t}X_{i}|\geq c_{k}\sqrt{n}+kt-1-t(k-1)>c_{k}\sqrt{|V(G')|}$.
Hence by Theorem~\ref{ThmC}, $G'$ has a $[2,k]$-ST $T$.
Then $T+\{u_{i}v:1\leq i\leq t,~v\in X_{i}\}$ is a spanning $(k,U)$-good tree of $G$.
\qed

\section{Proof of Theorems~\ref{thm2+} and \ref{thm3}}\label{sec-pf-thm3}

\begin{prop}
\label{prop-GGkn}
Let $k\geq 2$ and $n\geq 2k+1$ be integers.
Then for every $G\in \GG_{k,n}$, $G$ is a connected graph of order $n$ and satisfies $\sigma _{2}(G)=n-2$.
\end{prop}
\proof
Let $G\in \GG_{k,n}$.
By the definition of $\GG_{k,n}$, it is clear that $G$ is a connected graph of order $n$.
Let $L_{1}$, $L_{2}$, $L_{3}$ and $L_{4}$ be subsets of $V(G)$ satisfying (L1)--(L5).
Then the following hold.
\begin{enumerate}[{$\bullet $}]
\item
For every $u_{1}\in L_{1}$, it follows from (L2) and (L3) that $d_{G}(u_{1})=|L_{1}\cup L_{2}|-1$.
\item
For every $u_{2}\in L_{2}$, it follows from (L2) and (L4) that $d_{G}(u_{2})=|L_{1}\cup L_{2}|-1+|N_{G}(u_{2})\cap L_{3}|\geq |L_{1}\cup L_{2}|$.
Since $d_{G}(u_{2})\leq k$ for $u_{2}\in L_{2}$, this implies that $|L_{1}\cup L_{2}|\leq k$.
\item
For every $u_{3}\in L_{3}$, it follows from (L5) that $d_{G}(u_{3})\geq n-|L_{1}\cup L_{2}|-1$.
\item
For every $u_{4}\in L_{4}$, it follows from (L1)--(L5) that $d_{G}(u_{4})=|L_{3}\cup L_{4}|-1=n-|L_{1}\cup L_{2}|-1$.
\end{enumerate}
Since $n-|L_{1}\cup L_{2}|-1\geq (2k+1)-k-1=k\geq |L_{1}\cup L_{2}|-1$, it follows that
\begin{enumerate}[{$\bullet $}]
\item
for $u\in L_{1}\cup L_{2}$ and $u'\in L_{3}\cup L_{4}$ with $uu'\notin E(G)$, $d_{G}(u)+d_{G}(u')\geq (|L_{1}\cup L_{2}|-1)+(n-|L_{1}\cup L_{2}|-1)=n-2$,
\item
for $u_{1}\in L_{1}$ and $u_{4}\in L_{4}$, $d_{G}(u_{1})+d_{G}(u_{4})=(|L_{1}\cup L_{2}|-1)+(n-|L_{1}\cup L_{2}|-1)=n-2$, and
\item
for $u_{3},u'_{3}\in L_{3}$ with $u_{3}\neq u'_{3}$ and $u_{3}u'_{3}\notin E(G)$, $d_{G}(u_{3})+d_{G}(u'_{3})\geq 2(n-|L_{1}\cup L_{2}|-1)\geq (n-|L_{1}\cup L_{2}|-1)+(|L_{1}\cup L_{2}|-1)=n-2$.
\end{enumerate}
Consequently, we obtain $\sigma _{2}(G)=n-2$.
\qed

Now we prove Theorems~\ref{thm2+} and \ref{thm3}.

\medbreak\noindent\textit{Proof of Theorems~\ref{thm2+} and \ref{thm3}.}\quad
Let $k\geq 2$ be an integer.
For an integer $n\geq n_{0}(k)$, since $0\leq n-4c_{k}\sqrt{n}-2k^{2}-4k-4<n-2k$, i.e., $n\geq 2k+1$, it follows from Proposition~\ref{prop-GGkn} that each element of $\GG_{k,n}$ satisfies the assumption of Theorem~\ref{thm2+}.
Furthermore, if a graph $G$ satisfies (L1)--(L5), then $L_{2}$ is a cutset of $G$, and so $L_{2}$ is a $k$-blocking set of $G$ by (L4).
As we mentioned in Section~\ref{sec1}, if a graph has a $k$-blocking set, then the graph has no $[2,k]$-ST.
Therefore, to complete the proof of Theorems~\ref{thm2+} and \ref{thm3}, it suffices to show that a connected graph $G$ of order $n$ satisfying one of the following has a $[2,k]$-ST:
\begin{enumerate}[{\bf (G1)}]
\item[{\bf (G1)}]
$n\geq n_{0}(k)$, $\sigma _{2}(G)\geq n-2$ and $G$ is not isomorphic to any graph in $\GG_{k,n}$, or
\item[{\bf (G2)}]
(G1) does not hold, $n\geq n_{1}(k)$, $\sigma _{2}(G)\geq \frac{n+2k-2}{2}$ and $G$ has no $k$-blocking set.
\end{enumerate}

Since $k\geq 2$,
\begin{align}
c_{k}=\sqrt{k(k-1)(k+2\sqrt{2k}+2)}>\sqrt{k(k-1)(k+2)}\geq \sqrt{k^{3}}.\label{cond-ck}
\end{align}
By the definition of $n_{0}(k)$ and $n_{1}(k)$, we have
\begin{enumerate}[{\bf (N1)}]
\item[{\bf (N1)}]
if $n\geq n_{0}(k)$, then $n\geq 4c_{k}\sqrt{n}+2k^{2}+4k+4$, and
\item[{\bf (N'1)}]
if $n\geq n_{1}(k)$, then $\frac{n+2k-2}{4}\geq 2c_{k}\sqrt{n}+k^{2}+2k+1$.
\end{enumerate}
By (\ref{cond-ck}), if $k^{3}\geq n_{1}(k)$, then
\begin{align*}
0 &\leq \frac{k^{3}+2k-2}{4}-2c_{k}\sqrt{k^{3}}-k^{2}-2k-1\\
&< \frac{k^{3}+2k-2}{4}-2\sqrt{k^{3}}\cdot \sqrt{k^{3}}-k^{2}-2k-1\\
&= -\frac{7k^{3}+4k^{2}+6k+6}{4}\\
&< 0,
\end{align*}
which is a contradiction.
Thus $n_{1}(k)>k^{3}$.
In particular, if $n\geq n_{1}(k)$, then
$$
n-\frac{n+2k-2}{4}=\frac{3n-2k+2}{4}>\frac{3k^{3}-2k+2}{4}>0,
$$
and by (\ref{cond-ck}),
\begin{align}
c_{k}\sqrt{n}>\sqrt{k^{3}}\cdot \sqrt{k^{3}}=k^{3}\geq 2k^{2}.\label{eq-3-new1}
\end{align}
Consequently, it follows from (N'1) that
\begin{enumerate}[{\bf (N1)}]
\item[{\bf (N'2)}]
if $n\geq n_{1}(k)$, then $n>\frac{n+2k-2}{4}>c_{k}\sqrt{n}+3k^{2}+2k+1$.
\item[{\bf (N'3)}]
if $n\geq n_{1}(k)$, then $\frac{n+2k-2}{4}>5k^{2}+2k+1$.
\end{enumerate}

Let $s_{0}=n-2$ if (G1) holds; and let $s_{0}=\frac{n+2k-2}{2}$ if (G2) holds.
Let $S=\{u\in V(G):d_{G}(u)<\frac{s_{0}}{2}\}$.
Then for vertices $u,u'\in S$ with $u\neq u'$, we have $d_{G}(u)+d_{G}(u')<s_{0}\leq \sigma _{2}(G)$, and hence $uu'\in E(G)$.
This implies that $S$ is a clique of $G$.

If $\delta (G)\geq c_{k}\sqrt{n}$, then by Theorem~\ref{ThmC}, $G$ has a $[2,k]$-ST.
Thus we may assume that $\delta (G)<c_{k}\sqrt{n}$.
If (G1) holds, then by (N1), $\frac{s_{0}}{2}=\frac{n-2}{2}\geq \frac{(4c_{k}\sqrt{n}+2k^{2}+4k+4)-2}{2}>c_{k}\sqrt{n}$; if (G2) holds, by (N'1), $\frac{s_{0}}{2}=\frac{n+2k-2}{4}\geq 2c_{k}\sqrt{n}+k^{2}+2k+1>c_{k}\sqrt{n}$.
In either case, we have $\frac{s_{0}}{2}>c_{k}\sqrt{n}$, and hence there exists a vertex $u_{0}\in S$ such that $d_{G}(u_{0})=\delta (G)~(<c_{k}\sqrt{n})$.
In particular, $S\neq \emptyset $.
Since $S$ is a clique of $G$,
\begin{align}
|S|=|(N_{G}(u_{0})\cap S)\cup \{u_0\}|\leq d_{G}(u_{0})+1<c_{k}\sqrt{n}+1.\label{cond-3-1}
\end{align}

Let $\QQ$ be the family of components of $G-S$.
If $\QQ=\emptyset $, then by (\ref{cond-3-1}), $n=|S|<c_{k}\sqrt{n}+1$, which contradicts (N1) or (N'2).
Thus $\QQ\neq \emptyset $.
Furthermore, for each $Q\in \QQ $,
\begin{align}
\delta (Q)\geq \min\{d_{G}(v):v\in V(Q)\}-|S|\geq \frac{s_{0}}{2}-|S|=
\begin{cases}
\frac{n-2}{2}-|S| & \mbox{(if (G1) holds)}\\
\frac{n+2k-2}{4}-|S| & \mbox{(if (G2) holds)}.
\end{cases}
\label{cond-Q-s0}
\end{align}

\begin{claim}
\label{cl-3-0}
For $u\in S$ and $Q\in \QQ$, $V(Q)\setminus N_{G}(u)\neq \emptyset $.
\end{claim}
\proof
Suppose that $V(Q)\subseteq N_{G}(u)$.
Then $|V(Q)|+|S|-1=|V(Q)\cup (S\setminus \{u\})|\leq d_{G}(u)<\frac{s_{0}}{2}$.
On the other hand, for a vertex $v\in V(Q)$, we have $|V(Q)|+|S|-1=|(V(Q)\setminus \{v\})\cup S|\geq d_{G}(v)\geq \frac{s_{0}}{2}$, which is a contradiction.
\qed

Take $Q_{1}\in \QQ$ so that $|V(Q_{1})|$ is as small as possible.
Then
\begin{align}
|V(Q_{1})|\leq \frac{n-|S|}{|\QQ|}.\label{cond-3-new1}
\end{align}

\begin{claim}
\label{cl-3-1}
\begin{enumerate}[{\upshape(i)}]
\item
If (G1) holds, then $|\QQ|=1$.
\item
If (G2) holds, then $|\QQ|\leq 3$.
\item
If (G2) holds and $|S|\leq 3k-2$, then $|\QQ|\leq 2$.
\end{enumerate}
\end{claim}
\proof
Recall that $u_{0}$ is a vertex in $S$ with $d_{G}(u_{0})<c_{k}\sqrt{n}$.
By Claim~\ref{cl-3-0} with $(u,Q)=(u_{0},Q_{1})$, there exists a vertex $v_{1}\in V(Q_{1})\setminus N_{G}(u_{0})$.
In order to prove statements of the claim, we show that $d(u_0)+d(u_1)$ does not satisfy the degree-sum condition. 
By (\ref{cond-3-new1}),
\begin{align}
d_{G}(v_{1}) &\leq |V(Q_{1})\setminus \{v_{1}\}|+|S\setminus \{u_{0}\}|\nonumber \\
&\leq \left(\frac{n-|S|}{|\QQ|}-1\right)+(|S|-1)\nonumber \\
&= \frac{n+(|\QQ|-1)|S|}{|\QQ|}-2\label{cond-cl-3-1-1}.
\end{align}
To prove (i), (ii) and (iii), we prepare three equations as follows.
If $n\geq n_{0}(k)$, then by (N1),
\begin{align}
n-2-\left(c_{k}\sqrt{n}+\frac{n+c_{k}\sqrt{n}+1}{2}-2\right) &= \frac{n-3c_{k}\sqrt{n}-1}{2}\nonumber \\
&\geq \frac{(4c_{k}\sqrt{n}+2k^{2}+4k+4)-3c_{k}\sqrt{n}-1}{2}\nonumber \\
&\geq \frac{c_{k}\sqrt{n}+2k^{2}+4k+3}{2}\nonumber \\
&> 0.\label{cond-cl-3-1-2}
\end{align}
If $n\geq n_{1}(k)$, then by (N'1),
\begin{align}
\frac{n+2k-2}{2}&-\left(c_{k}\sqrt{n}+\frac{n+3c_{k}\sqrt{n}+3}{4}-2\right)\nonumber \\
&= \frac{n+2k-2}{4}-\frac{7c_{k}\sqrt{n}-2k-3}{4}\nonumber \\
&\geq (2c_{k}\sqrt{n}+k^{2}+2k+1)-\frac{7c_{k}\sqrt{n}-2k-3}{4}\nonumber \\
&= \frac{c_{k}\sqrt{n}+4k^{2}+10k+7}{4}\nonumber \\
&>0\label{cond-cl-3-1-3}
\end{align}
and
\begin{align}
\frac{n+2k-2}{2}&-\left(c_{k}\sqrt{n}+\frac{n+6k-4}{3}-2\right)\nonumber \\
&= \frac{2}{3}\cdot \frac{n+2k-2}{4}-\frac{3c_{k}\sqrt{n}+4k-8}{3}\nonumber \\
&\geq \frac{2(2c_{k}\sqrt{n}+k^{2}+2k+1)}{3}-\frac{3c_{k}\sqrt{n}+4k-8}{3}\nonumber \\
&= \frac{c_{k}\sqrt{n}+2k^{2}+10}{3}\nonumber \\
&> 0.\label{cond-cl-3-1-4}
\end{align}

Suppose that (G1) holds and $|\QQ|\geq 2$.
By (\ref{cond-3-1}), $\frac{n+(|\QQ|-1)|S|}{|\QQ|}\leq \frac{n+|S|}{2}<\frac{n+c_{k}\sqrt{n}+1}{2}$.
This together with (\ref{cond-cl-3-1-1}) and (\ref{cond-cl-3-1-2}) implies that
$$
\sigma _{2}(G)\leq d_{G}(u_{0})+d_{G}(v_{1})<c_{k}\sqrt{n}+\left(\frac{n+c_{k}\sqrt{n}+1}{2}-2\right)<n-2,
$$
which contradicts the assumption on $\sigma _{2}(G)$ in (G1).

Suppose that (G2) holds.
If $|\QQ|\geq 4$, then it follows from (\ref{cond-3-1}) that $\frac{n+(|\QQ|-1)|S|}{|\QQ|}\leq \frac{n+3|S|}{4}<\frac{n+3(c_{k}\sqrt{n}+1)}{4}$, and hence by (\ref{cond-cl-3-1-1}) and (\ref{cond-cl-3-1-3}),
$$
\sigma _{2}(G)\leq d_{G}(u_{0})+d_{G}(v_{1})<c_{k}\sqrt{n}+\left(\frac{n+3c_{k}\sqrt{n}+3}{4}-2\right)<\frac{n+2k-2}{2};
$$
if $|S|\leq 3k-2$ and $|\QQ|=3$, then $\frac{n+(|\QQ|-1)|S|}{|\QQ|}\leq \frac{n+2(3k-2)}{3}$, and hence by (\ref{cond-cl-3-1-1}) and (\ref{cond-cl-3-1-4}),
$$
\sigma _{2}(G)\leq d_{G}(u_{0})+d_{G}(v_{1})<c_{k}\sqrt{n}+\left(\frac{n+6k-4}{3}-2\right)<\frac{n+2k-2}{2}.
$$
In either case, we obtain a contradiction to the assumption on $\sigma _{2}(G)$ in (G2).
\qed

A set $S'\subseteq S$ {\it dominates} $\QQ$ if for each $Q\in \QQ$, there exists a vertex $u\in S$ such that $N_{G}(u)\cap V(Q)\neq \emptyset $.

\begin{claim}
\label{cl-3-case1-largedeg}
If there exists a vertex $u\in S$ such that $\{u\}$ dominates $\QQ$ and $d_{G}(u)\geq k+1$, then $G$ has a $[2,k]$-ST.
\end{claim}
\proof
Since $|N_{G}(u)\setminus (\bigcup _{Q\in \QQ}V(Q))|=|N_{G}(u)\cap S|=|S|-1$, there exists a set $X\subseteq N_{G}(u)\cap (\bigcup _{Q\in \QQ}V(Q))$ such that $|X|=\max\{|\QQ|,k+1-(|S|-1)\}$ and $X\cap V(Q)\neq \emptyset $ for all $Q\in \QQ$.
Note that $1\leq |X\cap V(Q)|\leq k+1$ for every $Q\in \QQ$.
Fix $Q\in \QQ$.
If (G1) holds, then by (N1), (\ref{cond-3-1}) and (\ref{cond-Q-s0}),
\begin{align*}
\delta (Q) &\geq \frac{n-2}{2}-|S|\\
&> \frac{n-2}{2}-c_{k}\sqrt{n}-1\\
&\geq \frac{(4c_{k}\sqrt{n}+2k^{2}+4k+4)-2}{2}-c_{k}\sqrt{n}-1\\
&= c_{k}\sqrt{n}+(k+1)^{2}-1\\
&> c_{k}\sqrt{|V(Q)|}+(k+1)|X\cap V(Q)|-1;
\end{align*}
if (G2) holds, then by (N'1), (\ref{cond-3-1}) and (\ref{cond-Q-s0}),
\begin{align*}
\delta (Q) &\geq \frac{n+2k-2}{4}-|S|\\
&> \frac{n+2k-2}{4}-c_{k}\sqrt{n}-1\\
&\geq (2c_{k}\sqrt{n}+k^{2}+2k+1)-c_{k}\sqrt{n}-1\\
&= c_{k}\sqrt{n}+(k+1)^{2}-1\\
&> c_{k}\sqrt{|V(Q)|}+(k+1)|X\cap V(Q)|-1.
\end{align*}
In either case, $\delta (Q)>c_{k}\sqrt{|V(Q)|}+(k+1)|X\cap V(Q)|-1$.
This together with Lemma~\ref{lem-3-2} with $(G,U)=(Q,X\cap V(Q))$ implies that there exists a spanning forest of $Q$ consisting of exactly $|X\cap V(Q)|$ components $F_{Q,1},F_{Q,2},\ldots ,F_{Q,|X\cap V(Q)|}$ such that for every integer $i$ with $1\leq i\leq |X\cap V(Q)|$, $|V(F_{Q,i})\cap X|=1$ and $F_{Q,i}$ is a $(k,V(F_{Q,i})\cap X)$-good tree.
Let
$$
T_{1}:=\left(\bigcup _{Q\in \QQ}\left(\bigcup _{1\leq i\leq |X\cap V(Q)|}F_{Q,i}\right)\right)+\{uv:v\in (S\setminus \{u\})\cup X\}.
$$
Then $d_{T_{1}}(u)=(|S|-1)+|X|\geq (|S|-1)+k+1-(|S|-1)=k+1$, and hence $T_{1}$ is a $[2,k]$-ST of $G$.
\qed

By Claim~\ref{cl-3-case1-largedeg}, we may assume that
\begin{align}
\mbox{if $\{u\}$ dominates $\QQ$, then $d_{G}(u)\leq k$.}\label{cond-one-dom-deg}
\end{align}

Note that $S$ dominates $\QQ$ because $G$ is connected.
Choose a set $\tilde{S}\subseteq S$ dominating $\QQ$ so that
\begin{enumerate}[{\bf (S1)}]
\item
$\tilde{S}$ is minimal, i.e., $\tilde{S}\setminus \{\tilde{u}\}$ does not dominate $\QQ$ for every $\tilde{u}\in \tilde{S}$,
\item
subject to (S1), $|\tilde{S}|$ is as large as possible, and
\item
subject to (S2), $\sum _{\tilde{u}\in \tilde{S}}d_{G}(\tilde{u})$ is as large as possible.
\end{enumerate}
If (G1) holds, then by Claim~\ref{cl-3-1}(i), $|\tilde{S}|=|\QQ|=1$; if (G2) holds, then by Claim~\ref{cl-3-1}(ii), $|\tilde{S}|\leq |\QQ|\leq 3$.
Write $\tilde{S}=\{\tilde{u}_{1},\tilde{u}_{2},\ldots ,\tilde{u}_{s}\}$ where $s=|\tilde{S}|$.

For $Q\in \QQ$, let $A_{Q}=\{u\in S:N_{G}(u)\cap V(Q)\neq \emptyset \}$.

\begin{claim}
\label{cl-3-2}
Let $Q\in \QQ$.
If (G2) holds and $\max\{d_{G}(u):u\in A_{Q}\}\leq k$, then $|\QQ|=1$ and $A_{Q}=S$.
\end{claim}
\proof
If $|\QQ|\geq 2$ or $A_{Q}\neq S$, then $A_{Q}$ is a cutset of $G$.
Since $G$ has no $k$-blocking set, this leads to the desired conclusion.
\qed

Now we divide the proof into two cases.

\medskip
\noindent
\textbf{Case 1:} $|\tilde{S}|=1$.

Note that $\tilde{S}=\{\tilde{u}_{1}\}$.
By (\ref{cond-one-dom-deg}), $d_{G}(\tilde{u}_{1})\leq k$.

\begin{claim}
\label{cl-3-case1-1}
We have $|\QQ|=1$, i.e., $\QQ=\{Q_{1}\}$.
\end{claim}
\proof
Suppose that $|\QQ|\geq 2$.
By Claim~\ref{cl-3-1}(i), we can assume (G2) holds.
By Claim~\ref{cl-3-2}, for each $Q\in \QQ$, there exists a vertex $u_{Q}\in A_{Q}$ with $d_{G}(u_{Q})\geq k+1$.
Note that $\{u_{Q}:Q\in \QQ\}$ dominates $\QQ$.
Take a minimal set $\tilde{S}_{0}\subseteq \{u_{Q}:Q\in \QQ\}$ dominating $\QQ$.
Since $|\tilde{S}|=1$, it follows from (S1) and (S2) that $|\tilde{S}_{0}|=1$.
However, $d_{G}(\tilde{u}_{1})\leq k$ and $d_{G}(\tilde{u})\geq k+1$ where $\tilde{u}$ is the unique element of $\tilde{S}_{0}$, which contradicts (S3).
Thus $|\QQ|=1$, and so $\QQ=\{Q_{1}\}$.
\qed

\begin{claim}
\label{cl-3-case1-2}
We have $A_{Q_{1}}=S$.
\end{claim}
\proof
It follows from (\ref{cond-one-dom-deg}) that
\begin{align}
\max\{d_{G}(u):u\in A_{Q_{1}}\}\leq k.\label{cond-cl-3-case1-2-1}
\end{align}
If (G2) holds, then by Claim~\ref{cl-3-2} and (\ref{cond-cl-3-case1-2-1}), $A_{Q_{1}}=S$.
Thus we may assume that (G1) holds.

Suppose that $A_{Q_{1}}\neq S$.
Let $L_{1}=S\setminus A_{Q_{1}}$, $L_{2}=A_{Q_{1}}$, $L_{3}=\bigcup _{u\in A_{Q_{1}}}(N_{G}(u)\cap V(Q_{1}))$ and $L_{4}=V(Q_{1})\setminus L_{3}$.
Note that $L_{1}$, $L_{2}$ and $L_{3}$ are non-empty sets.
By Claim~\ref{cl-3-case1-1}, $V(G)$ is the disjoint union of $L_{1}$, $L_{2}$, $L_{3}$ and $L_{4}$.

For $u_{2}\in L_{2}~(=A_{Q_{1}})$, it follows from (\ref{cond-cl-3-case1-2-1}) that
$$
|N_{G}(u_{2})\cap V(Q_{1})|\leq d_{G}(u_{2})-|L_{1}|\leq k-1
$$
and
\begin{align}
|L_{2}|\leq |L_{1}\cup L_{2}|-1\leq (|N_{G}(u_{2})\cup \{u_{2}\}|-1)-1\leq ((k+1)-1)-1.\label{cond-cl-3-case-1-2-2}
\end{align}
In particular, $|L_{3}|\leq \sum _{u_{2}\in L_{2}}|N_{G}(u_{2})\cap V(Q_{1})|\leq (k-1)^{2}$.
Hence by (N1) and (\ref{cond-cl-3-case-1-2-2}), $|L_{4}|=|V(Q_{1})\setminus L_{3}|=n-|L_{1}\cup L_{2}|-|L_{3}|\geq n-k-(k-1)^{2}>0$.
Hence $L_{4}\neq \emptyset $.

Let $u\in L_{1}$.
For $v_{3}\in L_{3}$, 
$$
n-2\leq \sigma _{2}(G)\leq d_{G}(u)+d_{G}(v_{3})=|(L_{1}\cup L_{2})\setminus \{u\}|+d_{G}(v_{3}),
$$
and hence $d_{G}(v_{3})\geq n-|L_{1}\cup L_{2}|-1$.
For $v_{4}\in L_{4}$,
\begin{align*}
n-2 \leq \sigma _{2}(G)&\leq d_{G}(u)+d_{G}(v_{4})\\
&\leq |(L_{1}\cup L_{2})\setminus \{u\}|+|V(Q_{1})\setminus \{v_{4}\}|\\
&= (|L_{1}\cup L_{2}|-1)+(n-|L_{1}\cup L_{2}|-1)\\
&= n-2,
\end{align*}
which forces $N_{G}(v_{4})=V(Q_{1})\setminus \{v_{4}\}$.
Since $v_{3}$ and $v_{4}$ are arbitrary, it follows from (\ref{cond-cl-3-case1-2-1}) that $(G,L_{1},L_{2},L_{3},L_{4})$ satisfies (L1)--(L5).
Consequently, $G$ is isomorphic to a graph in $\GG_{k,n}$, which contradicts (G1).
\qed

By Claims~\ref{cl-3-case1-1} and \ref{cl-3-case1-2}, for each $u\in S$, there exists a vertex $v_{u}\in N_{G}(u)\cap V(Q_{1})$.
Let $W=\{v_{u}:u\in S\}$.
Then
\begin{align}
|W|\leq |S|=|S\setminus \{\tilde{u}_{1}\}|+1=(d_{G}(\tilde{u}_{1})-|N_{G}(\tilde{u}_{1})\cap V(Q_{1})|)+1\leq k-1+1=k.\label{cond-case1-new1}
\end{align}
If (G1) holds, then by (N1), (\ref{cond-Q-s0}) and (\ref{cond-case1-new1}),
\begin{align*}
\delta (Q_{1}) &\geq \frac{n-2}{2}-|S|\\
&\geq \frac{n-2}{2}-k\\
&\geq \frac{(4c_{k}\sqrt{n}+2k^{2}+4k+4)-2}{2}-k\\
&> c_{k}\sqrt{n}+k^{2}-1\\
&> c_{k}\sqrt{|V(Q_{1})|}+k|W|-1;
\end{align*}
if (G2) holds, then by (N'2), (\ref{cond-Q-s0}) and (\ref{cond-case1-new1}),
\begin{align*}
\delta (Q_{1}) &\geq \frac{n+2k-2}{4}-|S|\\
&\geq \frac{n+2k-2}{4}-k\\
&> (c_{k}\sqrt{n}+3k^{2}+2k+1)-k\\
&> c_{k}\sqrt{n}+k^{2}-1\\
&> c_{k}\sqrt{|V(Q_{1})|}+k|W|-1.
\end{align*}
In either case, we have $\delta (Q_{1})>c_{k}\sqrt{|V(Q_{1})|}+k|W|-1$.
Hence by Lemma~\ref{lem-3-3} with $(G,U)=(Q_{1},W)$, there exists a spanning $(k,W)$-good tree $T$ of $Q_{1}$.
Then $T+\{uv_{u}:u\in S\}$ is a $[2,k]$-ST of $G$.

\medskip
\noindent
\textbf{Case 2:} $|\tilde{S}|\in \{2,3\}$.

By Claim~\ref{cl-3-1}(i), we can assume (G2) holds.
For the moment, we suppose that $|S|\geq (k-1)|\tilde{S}|+2$.
Recall that $s=|\tilde{S}|$.
Since $|S\setminus \tilde{S}|\geq (k-1)|\tilde{S}|-(|\tilde{S}|-2)$, there exists a partition $\{S_{1},S_{2},\ldots ,S_{s}\}$ of $S\setminus \tilde{S}$ such that
\begin{enumerate}[{$\bullet $}]
\item
if $|\tilde{S}|=2$, then $|S_{i}|\geq k-1$ for $i\in \{1,2\}$, and
\item
if $|\tilde{S}|=3$, then $|S_{i}|\geq k-1$ for $i\in \{1,3\}$ and $|S_{2}|\geq k-2$.
\end{enumerate}
Fix $Q\in \QQ$.
Since $\tilde{S}$ dominates $\QQ$, we can take an edge $u_{Q}v_{Q}\in E(G)$ with $u_{Q}\in \tilde{S}$ and $v_{Q}\in V(Q)$.
By (N'1), (\ref{cond-3-1}) and (\ref{cond-Q-s0}),
\begin{align*}
\delta (Q) &\geq \frac{n+2k-2}{4}-|S|\\
&> \frac{n+2k-2}{4}-c_{k}\sqrt{n}-1\\
&\geq (2c_{k}\sqrt{n}+k^{2}+2k+1)-c_{k}\sqrt{n}-1\\
&> c_{k}\sqrt{n}+k-1\\
&> c_{k}\sqrt{|V(Q)|}+k|\{v_{Q}\}|-1.
\end{align*}
Hence by Lemma~\ref{lem-3-3} with $(G,U)=(Q,\{v_{Q}\})$, there exists a spanning $(k,v_{Q})$-good tree $T_{Q}$ of $Q$.
Let $P=\tilde{u}_{1}\tilde{u}_{2}\cdots \tilde{u}_{s}$ be a path on $\tilde{S}$, and let
$$
T^{*}_{1}=\left(\left(\bigcup _{Q\in \QQ}T_{Q} \right)\cup P \right)+(\{u_{Q}v_{Q}:Q\in \QQ\}\cup \{\tilde{u}_{i}u:1\leq i\leq s,~u\in S_{i}\}).
$$
For $\tilde{u}_{i}\in \tilde{S}$, $|S_{i}|+d_{P}(\tilde{u}_{i})\geq k$ and, by (S1), there exists $Q\in \QQ$ such that $\tilde{u}_{i}=u_{Q}$, and hence $d_{T^{*}_{1}}(\tilde{u}_{i})=d_{P}(\tilde{u}_{i})+|\{Q\in \QQ:\tilde{u}_{i}=u_{Q}\}|+|S_{i}|\geq k+1$.
This implies that $T^{*}_{1}$ is a $[2,k]$-ST of $G$.
Thus we may assume that $|S|\leq (k-1)|\tilde{S}|+1~(\leq (k-1)|\QQ|+1)$.

If $|\QQ|=3$, then $|S|\leq 3k-2$, which contradicts Claim~\ref{cl-3-1}(iii).
Thus $|\QQ|=|\tilde{S}|=2$.
In particular,
\begin{align}
2=|\tilde{S}|\leq |S|\leq 2(k-1)+1=2k-1.\label{cond-case2-1}
\end{align}
Write $\QQ\setminus \{Q_{1}\}=\{Q_{2}\}$.
We may assume that $N_{G}(\tilde{u}_{i})\cap V(Q_{i})\neq \emptyset $ for each $i\in \{1,2\}$.
Then by (S1),
\begin{align}
\mbox{for $i\in \{1,2\}$, $N_{G}(\tilde{u}_{i})\cap V(Q_{3-i})=\emptyset $, i.e., $N_{G}(\tilde{u}_{i})\subseteq (S\setminus \{\tilde{u}_{i}\})\cup V(Q_{i})$.}\label{cond-case2-2}
\end{align}

\begin{claim}
\label{cl-3-3}
If $d_{G}(u)\leq 2k-1$ for every $u\in S$, then $G$ has a $[2,k]$-ST.
\end{claim}
\proof
It follows from (\ref{cond-case2-2}) that $N_{G}(\tilde{u}_{2})\cap V(Q_{1})=\emptyset$.
Hence by (\ref{cond-case2-1}) and the assumption of the claim, we have
\begin{align*}
\left|\bigcup _{u\in S}(N_{G}(u)\cap V(Q_{1}))\right| &\leq \sum _{u\in S\setminus \{\tilde{u}_{2}\}}|N_{G}(u)\cap V(Q_{1})|\\
&\leq \sum _{u\in S\setminus \{\tilde{u}_{2}\}}(d_{G}(u)-|S\setminus \{u\}|)\\
&\leq (|S|-1)(2k-1-(|S|-1))\\
&= (|S|-1)(2k-|S|)\\
&\leq (2k-2)^{2}.
\end{align*}
On the other hand, it follows from (N'3), (\ref{cond-Q-s0}) and (\ref{cond-case2-1}) that
\begin{align*}
|V(Q_{1})| &\geq \delta (Q_{1})+1\\
&\geq \frac{n+2k-2}{4}-|S|+1\\
&\geq \frac{n+2k-2}{4}-(2k-1)+1\\
&> (5k^{2}+2k+1)-(2k-1)+1\\
&> (2k-2)^{2}.
\end{align*}
Thus $V(Q_{1})\setminus (\bigcup _{u\in S}N_{G}(u))\neq \emptyset $.
Let $v^{*}\in V(Q_{1})\setminus (\bigcup _{u\in S}N_{G}(u))$.
Recall that we choose $Q_1$ so that $|V(Q_1)|$ is as small as possible.
Then by (\ref{cond-3-new1}), $d_{G}(v^{*})=|N_{G}(v^{*})\cap V(Q_{1})|\leq |V(Q_{1})\setminus \{v^{*}\}|\leq \frac{n-|S|}{2}-1$.

Fix $i\in \{1,2\}$.
Let $p_{i}=|N_{G}(\tilde{u}_{i})\cap V(Q_{i})|$.
Since $\tilde{u}_{i}v^{*}\notin E(G)$,
$$
\frac{n+2k-2}{2}\leq \sigma _{2}(G)\leq d_{G}(\tilde{u}_{i})+d_{G}(v^{*})\leq d_{G}(\tilde{u}_{i})+\frac{n-|S|}{2}-1,
$$
and hence $d_{G}(\tilde{u}_{i})\geq k+\frac{|S|}{2}$.
This together with (\ref{cond-case2-2}) implies that
\begin{align}
p_{i}=d_{G}(\tilde{u_{i}})-|S\setminus \{\tilde{u}_{i}\}|\geq k+\frac{|S|}{2}-(|S|-1)=k-\frac{|S|}{2}+1.\label{cond-cl-3-3-1}
\end{align}
Let $S'_{1}\subseteq S\setminus \tilde{S}$ be a set with $|S'_{1}|=\lceil \frac{|S\setminus \tilde{S}|}{2} \rceil $, and let $S'_{2}=S\setminus (\tilde{S}\cup S'_{1})$.
Note that $|S'_{1}|\geq |S'_{2}|=\lfloor \frac{|S\setminus \tilde{S}|}{2} \rfloor =\lfloor \frac{|S|-2}{2} \rfloor \geq \frac{|S|-3}{2}$.
Let $W_{i}=(N_{G}(\tilde{u}_{i})\cap V(Q_{i}))\cup S'_{i}$.
Then by (\ref{cond-cl-3-3-1}),
$$
|W_{i}|=p_{i}+|S'_{i}|\geq \left(k-\frac{|S|}{2}+1\right)+\frac{|S|-3}{2}=k-\frac{1}{2}.
$$
Since $|W_{i}|$ is an integer, this implies that $|W_{i}|\geq k$.
By (\ref{cond-case2-2}) and the assumption of the claim, we have
$$
|N_{G}(\tilde{u}_{i})\cap V(Q_{i})|=p_{i}=d_{G}(\tilde{u}_{i})-|N_{G}(\tilde{u}_{i})\cap S|\leq (2k-1)-|\{\tilde{u}_{3-i}\}|=2k-2.
$$
Hence by (N'2), (\ref{cond-Q-s0}) and (\ref{cond-case2-1}),
\begin{align*}
\delta (Q_{i}) &\geq \frac{n+2k-2}{4}-|S|\\
&\geq \frac{n+2k-2}{4}-(2k-1)\\
&> (c_{k}\sqrt{n}+3k^{2}+2k+1)-2k+1\\
&> c_{k}\sqrt{n}+(k+1)(2k-2)-1\\
&> c_{k}\sqrt{|V(Q_{i})|}+(k+1)|N_{G}(\tilde{u}_{i})\cap V(Q_{i})|-1.
\end{align*}
This together with Lemma~\ref{lem-3-2} with $(G,U)=(Q_{i},N_{G}(\tilde{u}_{i})\cap V(Q_{i}))$ implies that there exists a spanning forest of $Q_{i}$ consisting of exactly $p_{i}$ components $F'_{i,1},F'_{i,2},\ldots ,F'_{i,p_{i}}$ such that for every integer $j$ with $1\leq j\leq p_{i}$, $|V(F'_{i,j})\cap N_{G}(\tilde{u}_{i})|=1$ and $F'_{i,j}$ is a $(k,V(F'_{i,j})\cap N_{G}(\tilde{u}_{i}))$-good tree.
Then
$$
\left(\bigcup _{i\in \{1,2\}}\left(\bigcup _{1\leq j\leq p_{i}}F'_{i,j}\right)\right)+(\{\tilde{u}_{i}v:i\in \{1,2\},~v\in W_{i}\}\cup \{\tilde{u}_{1}\tilde{u}_{2}\})
$$
is a $[2,k]$-ST of $G$.
\qed

By Claim~\ref{cl-3-3}, we may assume that $\max\{d_{G}(u):u\in S\}\geq 2k$.
Since $|S|\leq 2k-1$, a vertex $u'\in S$ with $d_{G}(u')=\max\{d_{G}(u):u\in S\}$ satisfies $N_{G}(u')\setminus S\neq \emptyset $.
This together with (\ref{cond-one-dom-deg}) and (S3) forces $d_{G}(\tilde{u}_{i_{0}})\geq 2k$ for some $i_{0}\in \{1,2\}$.
Furthermore, it follows from Claim~\ref{cl-3-2} that $\max\{d_{G}(u):u\in A_{Q_{3-i_{0}}}\}\geq k+1$.
Hence by (\ref{cond-one-dom-deg}) and (S2), we have $d_{G}(\tilde{u}_{3-i_{0}})\geq k+1$.
Take a set $Z_{3-i_{0}}\subseteq N_{G}(\tilde{u}_{3-i_{0}})\setminus \{\tilde{u}_{i_{0}}\}$ such that $|Z_{3-i_{0}}|=k$ and $Z_{3-i_{0}}\cap V(Q_{3-i_{0}})\neq \emptyset $.
Then by (\ref{cond-case2-2}), $|N_{G}(\tilde{u}_{i_{0}})\cap Z_{3-i_{0}}|=|Z_{3-i_{0}}|-|Z_{3-i_{0}}\cap V(Q_{3-i_{0}})|\leq k-1$, and hence $|N_{G}(\tilde{u}_{i_{0}})\setminus (Z_{3-i_{0}}\cup \{\tilde{u}_{3-i_{0}}\})|\geq 2k-((k-1)+1)=k$.
In particular, we can take a set $Z_{i_{0}}\subseteq N_{G}(\tilde{u}_{i_{0}})\setminus (Z_{3-i_{0}}\cup \{\tilde{u}_{3-i_{0}}\})$ such that $1\leq |Z_{i_{0}}\cap V(Q_{i_{0}})|\leq k\leq |Z_{i_{0}}|$ and $S\setminus (Z_{3-i_{0}}\cup \tilde{S})\subseteq Z_{i_{0}}$.

Fix $i\in \{1,2\}$.
Then by the definition of $Z_{i}$, we have $q_{i}:=|Z_{i}\cap V(Q_{i})|\leq k$.
By (N'2), (\ref{cond-Q-s0}) and (\ref{cond-case2-1}),
\begin{align*}
\delta (Q_{i}) &\geq \frac{n+2k-2}{4}-|S|\\
&\geq \frac{n+2k-2}{4}-(2k-1)\\
&> (c_{k}\sqrt{n}+3k^{2}+2k+1)-2k+1\\
&> c_{k}\sqrt{n}+(k+1)k-1\\
&>c_{k}\sqrt{|V(Q_{i})|}+(k+1)|Z_{i}\cap V(Q_{i})|-1.
\end{align*}
Hence by Lemma~\ref{lem-3-2} with $(G,U)=(Q_{i},Z_{i}\cap V(Q_{i}))$, there exists a spanning forest of $Q_{i}$ consisting of exactly $q_{i}$ components $F_{i,1},F_{i,2},\ldots ,F_{i,q_{i}}$ such that for every integer $i$ with $1\leq j\leq q_{i}$, $|V(F_{i,j})\cap Z_{i}|=1$ and $F_{i,j}$ is a $(k,V(F_{i})\cap Z_{i})$-good tree.
Then
$$
\left(\bigcup _{i\in \{1,2\}}\left(\bigcup _{1\leq j\leq q_{i}}F_{i,j}\right)\right)+(\{\tilde{u}_{i}v:i\in \{1,2\},~v\in Z_{i}\}\cup \{\tilde{u}_{1}\tilde{u}_{2}\})
$$
is a $[2,k]$-ST of $G$.

This completes the proof of Theorem~\ref{thm3}.
\qed

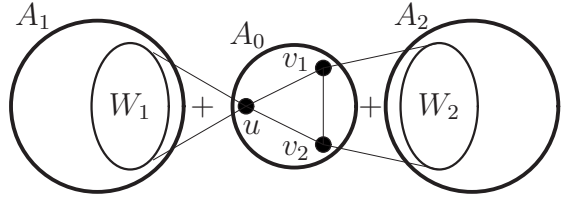
\begin{figure}
\begin{center}
{\unitlength 0.1in%
\begin{picture}(29.2200,9.8200)(0.9000,-12.4200)%
%
{\color[named]{Black}{%
\special{pn 0}%
\special{sh 1.000}%
\special{ia 1395 795 40 40 0.0000000 6.2831853}%
}}%
{\color[named]{Black}{%
\special{pn 8}%
\special{ar 1395 795 40 40 0.0000000 6.2831853}%
}}%
%
{\color[named]{Black}{%
\special{pn 0}%
\special{sh 1.000}%
\special{ia 1795 595 40 40 0.0000000 6.2831853}%
}}%
{\color[named]{Black}{%
\special{pn 8}%
\special{ar 1795 595 40 40 0.0000000 6.2831853}%
}}%
%
{\color[named]{Black}{%
\special{pn 0}%
\special{sh 1.000}%
\special{ia 1795 995 40 40 0.0000000 6.2831853}%
}}%
{\color[named]{Black}{%
\special{pn 8}%
\special{ar 1795 995 40 40 0.0000000 6.2831853}%
}}%
%
{\color[named]{Black}{%
\special{pn 20}%
\special{ar 2565 795 447 447 0.0000000 6.2831853}%
}}%
%
{\color[named]{Black}{%
\special{pn 20}%
\special{ar 625 795 447 447 0.0000000 6.2831853}%
}}%
%
{\color[named]{Black}{%
\special{pn 13}%
\special{ar 795 795 200 330 0.0000000 6.2831853}%
}}%
%
{\color[named]{Black}{%
\special{pn 13}%
\special{ar 2395 795 200 330 0.0000000 6.2831853}%
}}%
%
{\color[named]{Black}{%
\special{pn 4}%
\special{pa 1795 995}%
\special{pa 2395 1125}%
\special{fp}%
}}%
%
{\color[named]{Black}{%
\special{pn 4}%
\special{pa 2395 465}%
\special{pa 1795 595}%
\special{fp}%
}}%
%
{\color[named]{Black}{%
\special{pn 4}%
\special{pa 915 515}%
\special{pa 1395 795}%
\special{fp}%
}}%
%
{\color[named]{Black}{%
\special{pn 4}%
\special{pa 1395 795}%
\special{pa 915 1075}%
\special{fp}%
}}%
%
{\color[named]{Black}{%
\special{pn 4}%
\special{pa 1795 995}%
\special{pa 1795 595}%
\special{fp}%
\special{pa 1795 595}%
\special{pa 1395 795}%
\special{fp}%
}}%
%
{\color[named]{Black}{%
\special{pn 4}%
\special{pa 1395 795}%
\special{pa 1795 995}%
\special{fp}%
}}%
%
{\color[named]{Black}{%
\special{pn 20}%
\special{ar 1645 795 320 320 0.0000000 6.2831853}%
}}%
\put(7.9500,-7.9500){\makebox(0,0){{\color[named]{Black}{$W_{1}$}}}}%
\put(23.9500,-7.9500){\makebox(0,0){{\color[named]{Black}{$W_{2}$}}}}%
\put(2.9500,-3.2500){\makebox(0,0){{\color[named]{Black}{$A_{1}$}}}}%
\put(22.5500,-3.2500){\makebox(0,0){{\color[named]{Black}{$A_{2}$}}}}%
\put(14.0500,-4.2500){\makebox(0,0){{\color[named]{Black}{$A_{0}$}}}}%
\put(14.2500,-9.0500){\makebox(0,0){{\color[named]{Black}{$u$}}}}%
\put(16.5500,-5.6500){\makebox(0,0){{\color[named]{Black}{$v_{1}$}}}}%
\put(16.5500,-10.2500){\makebox(0,0){{\color[named]{Black}{$v_{2}$}}}}%
\put(11.7000,-8.0000){\makebox(0,0){{\color[named]{Black}{$+$}}}}%
\put(20.4000,-8.0000){\makebox(0,0){{\color[named]{Black}{$+$}}}}%
\end{picture}}%
\caption{Graph $G'_{k,n}$}
\label{f-Gkn}
\end{center}
\end{figure}

The degree-sum condition in Theorem~\ref{thm3} is best possible in a sense:
Let $k\geq 2$ be an integer, and let $n$ be an odd integer with $n\geq n_{1}(k)$.
Since $n\geq 2k^{2}+1$ by (N'2), we have $\frac{n-3}{2}-(k+1)=\frac{n-2k-5}{2}\geq \frac{(2k^{2}+1)-2k-5}{2}=(k+1)(k-2)\geq 0$, i.e., $\frac{n-3}{2}\geq k+1$.
Let $A_{0}$, $A_{1}$ and $A_{2}$ be vertex-disjoint complete graphs with $|V(A_{0})|=3$ and $|V(A_{1})|=|V(A_{2})|=\frac{n-3}{2}$.
Write $V(A_{0})=\{u,v_{1},v_{2}\}$, and for each $i\in \{1,2\}$, take a set $W_{i}\subseteq V(A_{i})$ with $|W_{i}|=k-1$.
Let $G'_{k,n}=(\bigcup _{0\leq i\leq 2}A_{i})+(\{uw:w\in W_{1}\}\cup \{v_{1}w,v_{2}w:w\in W_{2}\})$ (see Figure~\ref{f-Gkn}).
Then $G'_{k,n}$ is a connected graph of order $n$ and $\sigma _{2}(G'_{k,n})=\frac{n+2k-3}{2}$.
Furthermore, $G'_{k,n}$ has no $k$-blocking set.
Thus the following proposition gives a sharpness of Theorem~\ref{thm3}.

\begin{prop}
There exists no $[2,k]$-ST of $G'_{k,n}$.
\end{prop}
\proof
Suppose that $G'_{k,n}$ has a $[2,k]$-ST $T$.
Since $u$ is a cut-vertex of $G'_{k,n}$, we have $d_{T}(u)=k+1$; since $\{v_{1},v_{2}\}$ is a cutset of $G'_{k,n}$, we have $d_{T}(v_{i})=k+1$ for some $i\in \{1,2\}$.
This implies that $uv_{1}v_{2}u$ is a cycle of $T$, which contradicts the fact that $T$ is a tree.
\qed

\section{Proof of Theorem~\ref{thm4}}\label{sec-pf-thm4}

Since $p_{k}=\left(\frac{3c_{k}+\sqrt{c_{k}^{2}+4k^{2}+8k+4}}{2}\right)^{2}$, we have
\begin{align}
\sqrt{p_{k}}=\frac{3c_{k}+\sqrt{c_{k}^{2}+4k^{2}+8k+4}}{2}>2c_{k}.\label{eq-thm4-pk-1}
\end{align}
By (\ref{eq-thm4-pk-1}), we obtain the following:
\begin{enumerate}[{\bf (M1)}]
\item[{\bf (M1)}]
Since $c_{k}\geq 4$, we have $\sqrt{p_{k}}>c_{k}+\frac{1}{c_{k}}$, and hence
$$
c_{k}(\sqrt{p_{k}}-c_{k})-1>c_{k}\left(\left(c_{k}+\frac{1}{c_{k}}\right)-c_{k}\right)-1=0.
$$
\item[{\bf (M2)}]
Since $\sqrt{p_{k}}>\frac{c_{k}+\sqrt{c_{k}^{2}+4k}}{2}$, we have 
\begin{align*}
p_{k}-\frac{k\sqrt{p_{k}}}{\sqrt{p_{k}}-c_{k}} &= \frac{\sqrt{p_{k}}(\sqrt{p_{k}}(\sqrt{p_{k}}-c_{k})-k)}{\sqrt{p_{k}}-c_{k}}\\
&> \frac{\sqrt{p_{k}}\left(\frac{c_{k}+\sqrt{c_{k}^{2}+4k}}{2}\left(\frac{c_{k}+\sqrt{c_{k}^{2}+4k}}{2}-c_{k}\right)-k\right)}{\sqrt{p_{k}}-c_{k}}\\
&= 0.
\end{align*}
\item[{\bf (M3)}]
We have
\begin{align*}
&(\sqrt{p_{k}}-2c_{k})(\sqrt{p_{k}}-c_{k})\nonumber \\
&\quad = \left(\frac{3c_{k}+\sqrt{c_{k}^{2}+4k^{2}+8k+4}}{2}-2c_{k}\right)\left(\frac{3c_{k}+\sqrt{c_{k}^{2}+4k^{2}+8k+4}}{2}-c_{k}\right)\nonumber \\
&\quad = k^{2}+2k+1.
\end{align*}
\end{enumerate}

Let $S=\{u\in V(G):d_{G}(u)<\sqrt{p_{k}n}\}$.
Then for vertices $u,u'\in S$ with $u\neq u'$, we have $d_{G}(u)d_{G}(u')<p_{k}n\leq \pi _{2}(G)$, and hence $uu'\in E(G)$.
This implies that $S$ is a clique of $G$.

If $\delta (G)\geq c_{k}\sqrt{n}$, then it follows from Theorem~\ref{ThmC} that $G$ has a $[2,k]$-ST.
Thus we may assume that $\delta (G)<c_{k}\sqrt{n}$.
This together with (\ref{eq-thm4-pk-1}) implies that $\delta (G)<\sqrt{p_{k}n}$.
In particular, $S\neq \emptyset $.
Since $S$ is a clique of $G$, for a vertex $u_{0}\in S$ with $d_{G}(u_{0})=\delta (G)$, we have
\begin{align}
|S|=|(N_{G}(u_{0})\cap S)\cup \{u\}|\leq d_{G}(u_{0})+1<c_{k}\sqrt{n}+1.\label{cond-4-1}
\end{align}

Let $\QQ$ be the family of components of $G-S$.
If $\QQ=\emptyset $, i.e., $G=G[S]$, then $G$ is a complete graph of order at least $k+2$, and hence $G$ has a $[2,k]$-ST.
Thus we may assume that $\QQ\neq \emptyset $.
By (\ref{eq-thm4-pk-1}) and (\ref{cond-4-1}), for each $Q\in \QQ$,
\begin{align}
\delta (Q) \geq \min\{d_{G}(v):v\in V(Q)\}-|S|> \sqrt{p_{k}n}-(c_{k}\sqrt{n}+1).\label{cond-4-Q1}
\end{align}

\begin{claim}
\label{cl-4-0}
For $u\in S$ and $Q\in \QQ$, $V(Q)\setminus N_{G}(u)\neq \emptyset $.
\end{claim}
\proof
Suppose that $V(Q)\subseteq N_{G}(u)$.
Then $|V(Q)|+|S|-1=|V(Q)\cup (S\setminus \{u\})|\leq d_{G}(u)<\sqrt{p_{k}n}$.
On the other hand, for a vertex $v\in V(Q)$, we have $|V(Q)|+|S|-1=|(V(Q)\setminus \{v\})\cup S|\geq d_{G}(v)\geq \sqrt{p_{k}n}$, which is a contradiction.
\qed

Take $Q_{1}\in \QQ$ so that $|V(Q_{1})|$ is as small as possible.
Then
\begin{align}
|V(Q_{1})|\leq \frac{n-|S|}{|\QQ|}<\frac{n}{|\QQ|}.\label{cond-4-new1}
\end{align}

\begin{claim}
\label{cl-4-1}
We have $|\QQ|<\frac{\sqrt{n}}{\sqrt{p_{k}}-c_{k}}$, and in particular, $\sqrt{n}>\sqrt{p_{k}}-c_{k}$.
\end{claim}
\proof
By Claim~\ref{cl-4-0}, there exists a vertex $v\in V(Q_{1})$ such that $S\not\subseteq N_{G}(v)$, and hence
\begin{align}
|V(Q_{1})\cup S|\geq d_{G}(v)+2\geq \sqrt{p_{k}n}+2.\label{eq-cl-4-1-1}
\end{align}
If $|\QQ|\geq \frac{\sqrt{n}}{\sqrt{p_{k}}-c_{k}}$, then it follows from (\ref{cond-4-1}) and (\ref{cond-4-new1}) that
$$
|V(Q_{1})\cup S|<\frac{n}{|\QQ|}+|S|<\frac{n}{\frac{\sqrt{n}}{\sqrt{p_{k}}-c_{k}}}+(c_{k}\sqrt{n}+1)=\sqrt{p_{k}n}+1,
$$
which contradicts (\ref{eq-cl-4-1-1}).
\qed

A set $S'\subseteq S$ {\it dominates} $\QQ$ if for each $Q\in \QQ$, there exists a vertex $u\in S$ such that $N_{G}(u)\cap V(Q)\neq \emptyset $.
Note that $S$ dominates $\QQ$ because $G$ is connected.
Take a minimum set $\tilde{S}\subseteq S$ dominating $\QQ$, and write $\tilde{S}=\{\tilde{u}_{1},\tilde{u}_{2},\ldots ,\tilde{u}_{s}\}$ where $s=|\tilde{S}|$.
By Claim~\ref{cl-4-1} and the minimality of $\tilde{S}$,
\begin{align}
s\leq |\QQ|<\frac{\sqrt{n}}{\sqrt{p_{k}}-c_{k}}.\label{cond-4-sQ}
\end{align}

\begin{claim}
\label{cl-4-3}
For each $\tilde{u}\in \tilde{S}$, $d_{G}(\tilde{u})\geq sk+1$.
\end{claim}
\proof
Suppose that $d_{G}(\tilde{u})\leq sk$.
By Claim~\ref{cl-4-0}, there exists a vertex $v\in V(Q_{1})$ such that $v\tilde{u}\notin E(G)$.
By (\ref{cond-4-new1}) and (\ref{cond-4-sQ}), $d_{G}(v)<|V(Q_{1})\cup S|\leq \frac{n-|S|}{|\QQ|}+|S|\leq \frac{n-|S|}{s}+|S|$.
This together with (M1), (M2), (\ref{cond-4-1}) and (\ref{cond-4-sQ}) leads to
\begin{align*}
\pi _{2}(G) &\leq d_{G}(\tilde{u})d_{G}(v)\\
&< sk\left(\frac{n-|S|}{s}+|S|\right)\\
&= k(n+(s-1)|S|)\\
&< k\left( n+\left( \frac{\sqrt{n}}{\sqrt{p_{k}}-c_{k}}-1\right) (c_{k}\sqrt{n}+1)\right)\\
&= k\left(\left(1+\frac{c_{k}}{\sqrt{p_{k}}-c_{k}}\right)n-\frac{(c_{k}(\sqrt{p_{k}}-c_{k})-1)\sqrt{n}}{\sqrt{p_{k}}-c_{k}}-1\right)\\
&< \frac{k\sqrt{p_{k}}n}{\sqrt{p_{k}}-c_{k}}\\
&< p_{k}n,
\end{align*}
which is a contradiction.
\qed

By Claim~\ref{cl-4-3}, $|N_{G}(\tilde{u})\setminus \tilde{S}|\geq sk+1-(s-1)=sk-s+2$ for every $\tilde{u}\in \tilde{S}$.
Hence there exist $s$ disjoint subsets $W_{1},W_{2},\ldots ,W_{s}$ of $V(G)\setminus \tilde{S}$ such that for each integer $i$ with $1\leq i\leq s$, $W_{i}\subseteq N_{G}(\tilde{u}_{i})$ and
$$
|W_{i}|=
\begin{cases}
k+1 & \mbox{(if $s=1$)}\\
k & \mbox{(if $s\geq 2$ and $i\in \{1,s\}$)}\\
k-1 & \mbox{(if $s\geq 2$ and $2\leq i\leq s-1$)}.
\end{cases}
$$
Let $W=\bigcup _{1\leq i\leq s}W_{i}$.
Note that $|W|=sk-s+2$.

By the minimality of $\tilde{S}$, for each integer $i$ with $1\leq i\leq s$, there exists $D_{i}\in \QQ $ such that $(\bigcup _{v\in V(D_{i})}N_{G}(v))\cap \tilde{S}=\{\tilde{u}_{i}\}$.
We may assume that $V(D_{i})\cap W\neq \emptyset $.
For each $Q\in \QQ$, if $Q\in \{D_{i}:1\leq i\leq s\}$, then $1\leq |V(Q)\cap W|\leq k+1~(\leq k+(s-1)(k-2)+1)$; otherwise, $|V(Q)\cap W|\leq |W|-\sum _{1\leq i\leq s}|V(D_{i})\cap W|\leq sk-2s+2$.
In either case, we have $|V(Q)\cap W|\leq sk-2s+3$.
Let $\QQ_{1}=\{Q\in \QQ:V(Q)\cap W=\emptyset \}$.
For each $Q\in \QQ_{1}$, it follows from the definition of $\tilde{S}$, there exists an edge $\tilde{u}_{Q}v_{Q}~(\tilde{u}_{Q}\in \tilde{S},~v_{Q}\in V(Q))$ of $G$.
Let $W'=W\cup \{v_{Q}:Q\in \QQ_{1}\}$.

Fix $Q\in \QQ$.
Then $1\leq |W'\cap V(Q)|\leq sk-2s+3$, and hence by Claim~\ref{cl-4-1}, (M3), (\ref{cond-4-Q1}) and (\ref{cond-4-sQ}),
\begin{align*}
\delta (Q) &> (\sqrt{p_{k}}-c_{k})\sqrt{n}-1\\
&= c_{k}\sqrt{n}+(\sqrt{p_{k}}-2c_{k})\sqrt{n}-1\\
&> c_{k}\sqrt{|V(Q)|}+(\sqrt{p_{k}}-2c_{k})\sqrt{n}-1\\
&= c_{k}\sqrt{|V(Q)|}+\frac{(k+1)^{2}\sqrt{n}}{\sqrt{p_{n}}-c_{k}}-1\\
&= c_{k}\sqrt{|V(Q)|}+\frac{(k+1)((k-2)\sqrt{n}+3\sqrt{n})}{\sqrt{p_{n}}-c_{k}}-1\\
&> c_{k}\sqrt{|V(Q)|}+\frac{(k+1)((k-2)\sqrt{n}+3(\sqrt{p_{n}}-c_{k}))}{\sqrt{p_{n}}-c_{k}}-1\\
&> c_{k}\sqrt{|V(Q)|}+(k+1)(s(k-2)+3)-1\\
&\geq c_{k}\sqrt{|V(Q)|}+(k+1)|W'\cap V(Q)|-1.
\end{align*}
This together with Lemma~\ref{lem-3-2} with $(G,U)=(Q,W'\cap V(Q))$ implies that there exists a spanning forest of $Q$ consisting of exactly $|W'\cap V(Q)|$ components $F_{Q,1},F_{Q,2},\ldots ,F_{Q,|W'\cap V(Q)|}$ such that for every integer $i$ with $1\leq i\leq |W'\cap V(Q)|$, $|V(F_{Q,i})\cap W'|=1$ and $F_{Q,i}$ is a $(k,V(F_{Q,i})\cap W')$-good tree.
Let $H$ be the graph obtained from the path $\tilde{u}_{1}\tilde{u}_{2}\cdots \tilde{u}_{s}$ by joining $\tilde{u}_{1}$ and all vertices in $S\setminus (\tilde{S}\cup W')$.
Then
$$
\left(\left(\bigcup _{Q\in \QQ}\left(\bigcup _{1\leq i\leq |W'\cap V(Q)|}F_{Q,i}\right)\right)\cup H\right)+(\{\tilde{u}_{i}w:1\leq i\leq s,~w\in W_{i}\}\cup \{\tilde{u}_{Q}v_{Q}:Q\in \QQ_{1}\})
$$
is a $[2,k]$-ST of $G$.

This completes the proof of Theorem~\ref{thm4}.

\section*{Acknowledgment}

The authors would like to thank Professor Akira Saito for helpful comments.
This work was supported by the Research Institute for Mathematical Sciences, an International Joint Usage/Research Center located in Kyoto University.
This work is also supported by JSPS KAKENHI Grant numbers, 23K03204 (to M.F.), 	25K07107 (to S.T.) and research grant of Senshu University 2025 (to S.T.).

\end{document}